\documentclass [11pt]{amsart}
\usepackage{amscd}
\usepackage{amsxtra}
\usepackage{amsthm}
\usepackage{epsfig}
\usepackage{graphicx,color}
\usepackage{amssymb}
\usepackage[all]{xy}
\usepackage[ansinew]{inputenc}
\usepackage[T1]{fontenc}
\newtheorem{theorem}{\sc Theorem}
\newtheorem{proposition}{\sc Proposition}[section]
\newtheorem{lemma}{\sc Lemma}[section]
\newtheorem{corollary}{\sc Corollary}[section]
\theoremstyle{remark}
\newtheorem{definition}{\sc Definition}[section]

\newtheorem{example}{\it Example}[section]

\font\tmsb=msbm10 at12pt
\font\smsb=msbm7
\font\ssmsb=msbm5
\newfam\msbfam
\textfont\msbfam=\tmsb
\scriptfont\msbfam=\smsb
\scriptscriptfont\msbfam=\ssmsb
\def \1{\mathbb {1}}
\def \RM{\mathbb {R}}
\def \NM{\mathbb{N}}
\def \ZM{\mathbb{Z}}
\def \CM{\mathbb{C}}

\def \At {\mathcal A}
\def \Kt {\mathcal{K}}
\def \HM {\mathbb{H}}

\def \Bt {\mathcal{B}}
\def \Ct {\mathcal{C}}

\def \Et {\mathcal{E}}
\def \Ht {\mathcal{H}}
\def \Mt {\mathcal{M}}

\def \d{\partial}

\def\dt{\delta} 
\def\a{\alpha}
\def\b{\beta}
\def\e{\varepsilon}  
\def\g{\gamma}
\def\p{\varphi}

\def\l{\lambda}

\def \t{\tilde}

\def \to{\longrightarrow} 
\def \w{\wedge}

\def \Mk{\mathfrak{M}}

\def \< {{\langle }}
\def \> {{\rangle }}
\def \Ker {{\rm Ker}}

\newcommand{\Ft}{{\mathcal F}}

\newcommand{\OM}{{\mathcal O}}
\newcommand{\Ot}{{{\mathcal O} }}

\newcommand{\Lt}{{\mathcal L}}

\begin{document}
\title [Finiteness and constructibility in local analytic geometry]
{Finiteness and constructibility in\\ local analytic geometry}
\author[Mauricio D. Garay]{Mauricio D. Garay}
\date{Original version October 2006, Revised May 2007}
\address{SISSA/ISAS, via Beirut 4, 34014 Trieste,
Italy.}
\email{garay@sissa.it}


\parindent=0cm


\begin{abstract}{Using the Houzel finiteness theorem and the Whitney-Thom stratification theory, we show, in local analytic geometry,
that the direct image sheaves of relatively constructible sheaves have coherent direct images.}
\end{abstract}
\maketitle
\section*{Introduction}
In 1953, Cartan-Serre and Schwartz proved that the cohomology spaces of a coherent analytic sheaf on a compact complex analytic manifold
are finite dimensional (\cite{Cartan_Serre,Schwartz}).
In 1960, Grauert proved his direct image theorem stating that the direct images sheaves $R^kf_* \Ft$, associated to a coherent
analytic sheaf $\Ft$, are coherent provided that the holomorphic map $f:X \to S$ is proper. It is only with the work of
Forster-Knorr and of Kiehl-Verdier, that a proof similar to the absolute one was done \cite{Forster_Knorr,Verdier}
(see also \cite{Douady,Levy}).
Their proof was simplified and generalised to more general sheaves by Houzel (\cite{Houzel}).
The aim of this paper is to deduce from the Houzel theorem, a practical criterion for the coherence of direct image sheaves, close in spirit to
the works \cite{Br3,Buchweitz,Houzel_Schapira,Schneiders,vanStraten}. The main difference of our approach is that, because we have in
view applications to singularity theory, we use
the Whitney-Thom theory of stratified sets and mapping, which will consequently occupy most of the paper.
One of the key results is the following theorem.
\begin{theorem}
\label{T::Grauert}
 {Let $f:(\CM^n,0) \to (\CM^k,0)$ be a holomorphic map germ satisfying the $a_f$-condition.
The cohomology spaces $H^p(K^\cdot) $ associated to a complex of  $f$-constructible
$\Ot_{\CM^n,0}$-coherent modules are $f^{-1}\Ot_{\CM^k,0}$-coherent modules, for any $p \geq 0$.
Moreover, the germ $f$ admits representatives
$f:X \to S$ such that the direct image sheaves are $\RM^\cdot f_* \Kt^\cdot$ are $\Ot_S$-coherent and such
that $H^p(K^\cdot)=(\RM^p f_* \Kt^\cdot)_0 $, for any $p \geq 0$.}
\end{theorem}
In the statement of the theorem $\Kt^\cdot$ is the complex of coherent analytic sheaves in $X$ whose stalk at the origin is equal to $K^\cdot$.\\
Here $f$-constructible means fibrewise constructible, a notion that we shall carefully explain in the sequel. The representatives of the theorem are called {\em standard representatives}, they will be defined in Section \ref{SS::standard}.\\
A particular case is when $f$ defines an isolated singularity and the complex is the relative de Rham complex. In particular for hypersurface singularities, i.e. for $k=1$,
we get the Brieskorn coherence theorem (\cite{Br3}). The proof of the theorem is indeed similar to that of Brieskorn.\\
Then, we will give a generalisation
of this theorem which is applicable for instance to complexes of modules over rings of pseudo-differential operators. As an illustration, we will give a simple proof of the algebraic form of the Sato-Kashiwara-Kawai division theorem for
pseudo-differential operators due to Boutet de Monvel.\\
The results of this paper might be well known to some specialists, but we have thought that a paper giving an elementary
presentation of the subject together with simple criteria, based on stratification theory, of the abstract theorems 
might be of some use.

\section{The finiteness theorem in the absolute case}
\subsection{Statement of the theorem}
Given a sheaf $\Ft$ in $\CM^n$, we denote by $\Ft_0$ its stalk at the origin.
We denote by $B_r \subset \CM^n$ the closed ball of radius $r$
centred at the origin and by $\dot B_r$ its interior.
In the absolute case Theorem \ref{T::Grauert} can be stated as follows. 
\begin{theorem}
\label{T::absolute}
 {For any  constructible complex $\Kt^\cdot$ of coherent
analytic sheaves in $B_r \subset \CM^n$, the cohomology spaces $H^p(K^\cdot), K^\cdot=\Kt_0^\cdot $,
are finite dimensional vector spaces, for any $p \geq 0$. Moreover, for $\e<r$ small enough, the canonical mapping
$\Kt^\cdot(B_\e) \to \Kt^\cdot_0=K$ induces an isomorphism
$$H^p(K^\cdot) \approx \HM^p(B_\e,\Kt^\cdot), \forall p \geq 0 $$}
\end{theorem}

Here the constructibility of the complex $\Kt^\cdot$ means that the cohomology sheaves $\Ht^k(\Kt^\cdot)$ are locally constant on the stratum
of some Whitney stratification.
The proof of this theorem is a simple variant of the original Cartan-Serre-Schwartz proof. Although it is quite elementary, it
contains in essence all the ingredients involved in the proofs of more elaborated results.
We first give an example of application.
\subsection{Finiteness of de Rham cohomology of an isolated singularity}
Consider the complex $\Omega^\cdot_{X}$ of K\"ahler differentials on a Stein complex variety $X \subset \CM^n$.
For instance, if $X$ is a hypersurface
then the terms of the complex are $\Omega^k_{X}=\Omega^k_{\CM^n}/(df \w \Omega^{k-1}_{\CM^n}+f\Omega^k_{\CM^n}) $,
where $f$ is a generator of the ideal of $X$, the differential of the complex
being induced by the de Rham differential.\\
The Poincar\'e lemma states that at the smooth points of $X$ the complex is a resolution of the constant sheaf $\CM_X$, therefore if $X$ has isolated singular points the complex is constructible. Applying Theorem \ref{T::absolute}, we get the following result.
\begin{proposition} If $(X,0) \subset (\CM^n,0)$ is the germ of a variety with an isolated singular point at the origin
then the complex of K\"ahler differentials has finite dimensional cohomology spaces. 
\end{proposition}
If $X$ is a hypersurface then, as conjectured by Brieskorn, these cohomology spaces are all zero except, possibly, for $j=0,n-1,n$ (\cite{Sebastiani}).
\subsection{A small review on the Whitney Stratification theory}
\label{SS::standard}
A pair of $C^\infty$ submanifolds $U,V \subset \RM^k$, $\dim V<\dim U$,  satisfies the {\em Whitney condition}
if the following property holds: for any pair of sequences $(x_i),(y_i)$ of the submanifolds $U$ and $V$
both converging to the same point, such that:
\begin{enumerate}
\item the sequences of secants $(x_iy_i)$ converges to a line $L$,
\item the sequence of spaces tangent to $U$ at $x_i$ converge to an affine subspace $A \subset \RM^k$
\end{enumerate}
then the line $L$ is contained in the affine subspace $A$.
\begin{definition}
A stratification $\bigcup_{i=1}^m X_i$ of a subset $X \subset \RM^k$ by $C^\infty$ manifolds is called a {\em Whitney stratification} if for any
stratum $X_j$ lying on the closure of a stratum $X_i$
the pair $(X_i,X_j)$ satisfies the Whitney condition.  
\end{definition}
These definitions are due to Whitney  \cite{Wh} (see also
\cite{Gibson_Looijenga}, Chapter 1). Whitney proved the existence of such a stratification for
any real semi-analytic subset, constructive proofs were given in \cite{Lojasiewicz,Te2}. Topological characterisation of Whitney
stratification were given in \cite{Teissier_Oslo}.\\
We shall say that two Whitney stratified sets intersect {\em transversally} if their strata intersect  pairwise transversally.
We denote by $B_\e \subset \RM^k$ the closed ball centred at the origin of radius $\e$.
A direct consequence of the definition is the following.
\begin{proposition}
\label{P::transversalite}
Let $X \subset \RM^k$ be a Whitney stratified subset, then there exists $\e_0$ such that
the balls $B_\e, \e<\e_0$, intersect $X$ transversally.
\end{proposition}
\begin{proof}
If such an $\e_0$ does not exists then we can construct a sequence $(x_i)$ lying on a stratum, such that the affine space $T_i$
tangent to the stratum of $x_i$ at the point $x_i$ is also tangent to the boundary of the ball $B_{1/i}$ of radius $1/i$.
In particular the secant
$(Ox_i)$ is perpendicular to $T_i$. This contradicts the Whitney condition. 
\end{proof}
We say that a real vector field on a Whitney stratified set is {\em stratified} if its flow is continuous and if the restriction of the vector field to each stratum is $C^\infty$. A typical example, in $\RM^3$, of a stratified vector field
is given by 
$$\left\{  \begin{matrix} \d_z+\frac{x\d_y-y\d_x}{\sqrt{x^2+y^2}} & {\rm if ~ } (x,y) \neq (0,0) \\
                             \d_z &{\rm otherwise}
 \end{matrix} \right.$$
 for the stratification consisting of the $z$-axis and its complement.\\
A vector field will be called {\em tangent} to a stratified variety if it is tangent to all its strata.
One of the main property of Whitney stratifications is the following. 
\begin{theorem}[\cite{Wh}]
\label{T::Whitney}
 Any $C^\infty$ vector field tangent to the stratum $M$ of a point $x \in X$ of a Whitney stratified set$X$
extends to a stratified vector field tangent to the strata to which $M$ is adjacent.
\end{theorem}
This theorem can also be deduced from
\cite{Gibson_Looijenga}, Chapter II, Corollary 2.7 (see also Subsection \ref{SS::Thom}).
It implies in particular the germ of $X$ at a point which belongs to a stratum of positive dimension is homeomorphic to the germ of the product of a stratified manifold by a real line.
\begin{example}Consider the real singular surface
$$S=\{ (x,y,z) \in \RM^3: xy(x+y)(x-zy)=0 \}.$$
According to the Whitney theorem, in a small neighbourhood of the point $(x,y,z)=(0,0,1)$, it is
homeomorphic to the product of four lines by $\RM$. Different slices $\{z={\rm constant} \}$
consists of four lines having different cross-ratios, thus the surface is not locally diffeomorphic to a product.
\end{example}
\subsection{Construction of the contraction}
We apply the previous consideration to $\CM^n \approx \RM^{2n}$ with the stratification given by a constructible complex of coherent
analytic sheaves $\Kt^\cdot$ defined in a neighbourhood $U \subset \CM^n$ of the origin. Recall that the complex is {\em constructible} means that there exists a Whitney stratification of the neighbourhood $U$ such that the cohomology
sheaves are locally constant on each stratum. According to Proposition \ref{P::transversalite}, we can find a ball
$B_{\e_0} \subset U$ such that all strata in $U$ are transverse to the boundaries of the balls
$B_{\e}$ for any $\e<\e_0$.\\
The aim of this subsection is to prove the following proposition.
\begin{proposition}
\label{P::contraction}
The restriction mappings $r:\Kt^\cdot( B_{\e}) \to \Kt^\cdot( B_{\e'})$,\\
$r':\Kt^\cdot(B_{\e}) \to \Kt^\cdot(\dot B_{\e})$ are quasi-isomorphisms for any $\e'<\e<\e_0$.
\end{proposition}
We consider only the case of $r$, the other case being perfectly similar.\\
It is obviously sufficient to prove the proposition for $\e'$ sufficiently close to $\e$.
\begin{lemma}For any $\e < \e_0$, there exists a stratified vector field in $B_{\e_0} \subset \CM^n$ 
everywhere transversal to the boundary of $ B_\e$.
\end{lemma}
\begin{proof}
Let $M \subset B_{\e_0}$ be the stratum of minimal dimension which intersects the ball $B_\e$. Let $d:B_{\e_0} \to \RM$
be the algebraic euclidean distance of a point $y \in B_{\e_0}$ to $B_\e$, that is, inside $B_\e$ it is the euclidean distance and outside $B_\e$, it is minus the euclidean distance.\\
We denote by $\theta_M$ the gradient vector field of $d$ with respect to the metric induced by $\CM^n \approx \RM^{2n}$
on $B_{\e_0}$. The vector field $\theta_M$ is a $C^\infty$ vector field tangent to $M$ and transversal to the intersection of $M$
with $B_\e$.\\
By Theorem \ref{T::Whitney}, 
the vector field $\theta_M$ extends to a stratified vector field $\theta$
in $B_{\e_0}$ transversal to $M \cap B_{\e_0}$.
\end{proof}

\begin{lemma}
\label{L::qisomorphism}
 The restriction mapping $r:\Kt^\cdot(B_\e) \to \Kt^\cdot(B_{\e'})$ is a quasi-isomorphism.
\end{lemma}
\begin{proof}
Denote by $\p:]-\dt,\dt[ \times B_\e \to B_\e$ the flow of the vector field $\theta$ where $\dt$ is small enough so that it
induces a map which is a homeomorphism onto a neighbourhood of the boundary of $B_\e$ as stratified sets:
$$]-\dt,\dt[ \times \d B_{\e} \to B_\e,\ (t,x) \mapsto \p(t,x).$$ 
We chose $\e'$ close enough to $B_\e$ so that the boundary of $B_{\e'}$ is contained in the image of this map.\\
Chose an acyclic covering $U=(U_i)$ of $B_\e$, its image $U'=(U_i')$, $U_i'=U_i \cap B_{\e'}$, is an acyclic covering of $B_{\e'}$.\\
Consider the spectral sequences $E_0^{p,q}(B_\e)=\Ct^p(U,\Kt^q)$,
$E_0^{p,q}(B_{\e'})=\Ct^p(U',\Kt^q)$ for the hypercohomology
of the complex $\Kt^\cdot$. Here, as usual, $\Ct^\cdot(\cdot)$ stands for the
\v{C}ech resolution.\\
 The map $\p$ induces a homeomorphism between each stratum in $U_i$ and the
corresponding stratum in $U_i'$ for each $i$.
As the cohomology sheaves are locally constant on each stratum, we have a group isomorphism $\Ht^q(\Kt^\cdot)(U_i) \approx \Ht^q(\Kt^\cdot)(U_i')$ on each small open subset $U_i$.
Therefore, the restriction mapping induces an isomorphism between the $E_1$-terms of the hypercohomology
spectral sequences:
$$E_1^{p,q}(B_\e)=\Ct^p(U,\Ht^q(\Kt^\cdot)) \approx \Ct^p(U',\Ht^q(\Kt^\cdot))=E_1^{p,q}(B_{\e'}).$$
This shows that the hypercohomology spaces
$\HM^\cdot(B_\e,\Kt^\cdot)$ and $\HM^\cdot(B_{\e'},\Kt^\cdot)$ are isomorphic.\\
By Cartan's theorem B, for any $p \geq 0$, we have the isomorphisms
$$\HM^p(B_\e,\Kt^\cdot) \approx H^p(\Kt^\cdot(B_\e)),\ \HM^p(B_{\e'},\Kt^\cdot) \approx
H^p(\Kt^\cdot(B_{\e'})),$$
therefore the restriction mapping is a quasi-isomorphism. This proves the lemma.
\end{proof}

\subsection{A small review on Grothendieck's theory of topological tensor products}
The restriction mapping constructed in the preceding subsection has the property to be nuclear, a notion due to Grothendieck that
we shall now explain.
We consider only vector spaces over the field of complex numbers.\\
A topological vector space is called {\em locally convex} if its topology is generated by a set of continuous semi-norms
$(p_n),\ n \in \Omega$, that is, the subsets $V_{n,\e}=\{ x \in E: p_n(x) <\e \}$ form a fundamental system of $0$-neighbourhoods.\\
A locally convex topological vector space $E$ is called  a {\em Fr\'echet space}  (or an $F$-space) if it is complete and metrisable or equivalently when the topology of $E$ is generated by a countable set of semi-norms.
\begin{example}
Consider the vector space $\Ot_{\CM}(D)$ of holomorphic functions on the open disk.
Each compact neighbourhood $K \subset U$ defines a norm $\| f \|_K=\sup_{x \in K} |f(x)|$.
The neighbourhoods $V_{K,\e}=\{ f; \| f \|_K<\e \}$ form a basis of neighbourhood of the origin which
induces a Fr\'echet space structure on the vector space $\Ot_\CM(U)$. Indeed, consider
the sequence $(K_n)$ of closed disks of radius $1-1/n$. The neighbourhoods $V_{K_n,\e}$ form a countable
basis of $0$-neighbourhoods in $\Ot_{\CM}(D)$. There is no difficulty in extending this example to higher dimensions.
\end{example}
A subset of a locally convex topological vector space is called {\em bounded} if all semi-norms are bounded on it.\\
The {\em strong dual} $E'$ of a topological vector space $E$ means the topological dual together with the topology
induced by the semi-norms
$$p_{B}(u)=\sup_{x \in B}|u(x)|$$
where $B$ runs over the bounded subsets of $E$. For instance if $E$ is Banach, this is the topological dual with the operator norm topology. These definitions are standard \cite{Bourbaki_evt}.\\
Given two locally convex spaces $E,F$, we denote by $E \hat \otimes F$
the set of expression of the type  $v=\sum \l_i x_i \otimes y_i$ where
the sequences $(x_i)$ and $(y_i)$ are bounded and $\sum_i | \l_i |<\infty$, it is called the {\em topological tensor product} of the two
spaces. 
\begin{theorem}[\cite{Grothendieck_these,Grothendieck_PTT}] The topological tensor product of two Fr\'echet spaces is complete and separated,
thus it is also
a Fr\'echet space.
\end{theorem}
\begin{theorem}[\cite{Grothendieck_these,Grothendieck_PTT}]
 The topological tensor product of the strong dual of a Fr\'echet space with a Fr\'echet space
is complete and separated.
\end{theorem}
A bounded linear mapping $u:E \to F$ between Fr\'echet spaces is called {\em nuclear}
if it lies in the image of the canonical (Kronecker) embedding $E' \hat \otimes F \to L(E,F)$. Nuclear mappings are limits of finite
range mappings, they are therefore compact.
\begin{example}
\label{E::nuclear}
 Take $E=\Ot_{\CM}(D)$, $F=\Ot_{\CM}(D')$ where
$D,D'$ are open disks centred at the origin such that the radius $r$ of the disk $D' \subset \CM$ is strictly less
than that the radius $R$ of $D$.\\
The restriction mapping $\rho:\OM_{\CM}(D) \to \OM_{\CM}(D')$ is nuclear.
To see it, define the linear forms $a_n$ 
$$a_n:\Ot_{\CM}(D) \to \CM,\ f \mapsto \frac{1}{2i\pi}\int_\g \frac{f(z)}{z^{n+1}} dz $$
where $\g$ is a path in $D \setminus D'$ which turns counterclockwise around the origin.
We have the equality
$$\rho=\sum_{n \geq 0} \l^n a_n \otimes \frac{z^n}{\l^n},\ \l=\frac{r}{R}.$$ 
This shows that $\rho$ is nuclear.
This result extends to arbitrary Stein neighbourhoods \cite{Verdier} (see also \cite{Hubbard}).
\end{example}
Consider a coherent analytic sheaf $\Ft$ in $\CM^n$ and take a presentation of this sheaf
$$\Ot^{n_1}_{\CM^n} \to \Ot^{n_0}_{\CM^n} \to \Ft \to 0$$
This exact sequence induces a Fr\'echet structure on the vector space $\Ft(U)$ for any Stein neighbourhood $U \subset \CM^n$.
This structure is independent on the choice of the presentation
\cite{Cartan_Serre} (see also \cite{DouadyR.}, Proposition 4).\\
The proof of the following proposition is a generalisation of our previous example.
\begin{proposition}
[\cite{Cartan_Serre,Verdier}]
\label{P::nuclear}
For any coherent analytic sheaf $\Ft$ in $\CM^n$ and any Stein neighbourhoods $U,U'$ so that the closure of $U'$ is a compact
subset of $U$, the restriction mapping $\Ft(U) \to \Ft(U')$ is nuclear.
\end{proposition}
\subsection{The Schwartz perturbation theorem and the Houzel lemma}
The following theorem was proved by Schwartz for the more general case of compact operators.
 \begin{theorem}[\cite{Schwartz}]
\label{T::Schwartz_perturbation}Let  $f:E \to F$ be a surjective continuous linear morphism between Fr\'echet spaces.
For any nuclear map $u:E \to F$, the cokernel of the map $f+u$ is of finite dimension.
\end{theorem}

\begin{corollary}[\cite{Schwartz}]
\label{C::Schwartz}Let  $M^\cdot, N^\cdot$ two complexes of Fr\'echet spaces
for which there exists a nuclear quasi-isomorphism $u:M^\cdot \to N^\cdot$.
Then, the cohomology spaces of the complexes are finite dimensional.
\end{corollary}
\begin{proof}
We apply Theorem \ref{T::Schwartz_perturbation} to the maps
$$f^k:M^k \times N^{k-1} \to N^k, (\a,\b) \mapsto u(\a)+d\b.$$
\end{proof}

The following result was proved by Houzel under much more general assumptions, for a proof of this lemma 
we refer to Houzel's paper.
\begin{lemma}[\cite{Houzel}]
\label{L::Houzel}
 Any nuclear  mapping $u:E \to E$ of a Fr\'echet space to itself can be written as $u=u'+u''$
where $u'$ is of finite range,
i.e. $\dim {\rm Im}\, u'<+\infty$, and $I+u''$ is invertible.
\end{lemma}
Remark that for Banach spaces the proof of the lemma is obvious.
This lemma implies the Schwartz perturbation theorem. Indeed, consider first the case where $f$ is the identity mapping.
Take $u,u',u''$ as in the lemma, as the map $I+u''$ is invertible,
the map $I+u=(I+u'')+u'$
has a cokernel of finite dimension.\\
Now, consider an arbitrary surjective map $f:E \to F$.
I assert that any nuclear mapping $u$ factors through $f$ by a nuclear mapping:
$$ \xymatrix{ E \ar[rd]^-u \ar@{.>}[d]_v \\
E \ar[r]^f & F }
$$
As $u$ is nuclear, we may write $u=\sum_{i >0}\l_i \xi_i \otimes y_i$ where $\xi_i$ and $y_i$ are bounded sequences,
 $\sum_{i >0}|\l_i|<+\infty$ and $\displaystyle{\lim_{i \to +\infty} y_i=0}$.\\
As $f$ is surjective, the Banach open mapping theorem implies that there exists a sequence $(x_i)$ converging to zero in $E$
which lifts the sequence $(y_i)$, that is, $f(x_i)=y_i$ and
$\displaystyle{\lim_{i \to +\infty} x_i=0}$.\\
We define the nuclear map $v$ by the formula $v=\sum_{i >0}\l_i \xi_i \otimes x_i$. The map $I+v$ and consequently the
map $f \circ (I+v)=f+u$ have a finite dimensional cokernel. This shows that the Schwartz perturbation theorem is a consequence
of Lemma \ref{L::Houzel}.
\subsection{Proof of Theorem \ref{T::absolute}}

As the restriction mapping $\Kt^\cdot(\dot B_\e) \to \Kt^\cdot(\dot B_{\e'})$ is a nuclear quasi-isomorphism (Proposition \ref{P::contraction} and Proposition \ref{P::nuclear}),
Corollary \ref{C::Schwartz} applies.
This shows that the cohomology spaces of the complex $\Kt^\cdot(\dot B_\e)$ are finite dimensional vector spaces or
equivalently that it is quasi-isomorphic to a complex $\Lt^\cdot$ of finite dimensional constant sheaves $\Lt^k \approx \CM^{n_k}$. \\
As  $\Kt^\cdot(\dot B_\e)$ is quasi-isomorphic to $\Kt^\cdot(\dot B_{\e'})$,
we may construct an exhaustive sequence of compact neighbourhoods of the origin
$(B_{\e_n})$ so that $\Kt^\cdot(\dot B_{\e_n})$ is quasi-isomorphic to $\Kt^\cdot(\dot B_{\e_{n+1}})$.
In the limit $n \to \infty$, we get that the complex $K^\cdot=\Kt^\cdot_0$ is quasi-isomorphic to the stalk of the complex $\Lt^\cdot$ at the origin.
This concludes the proof of the theorem in the absolute case.

\section{Relatively constructible sheaves}
We explain the notion of $f$-constructibility introduced in Theorem \ref{T::Grauert}.
\subsection{A review on stratified mappings}
\label{SS::Thom}
\begin{definition}A continuous map between Whitney stratified topological spaces $f:X \to Y$, $X=\cup_jX_j$
, $Y=\cup_j Y_j$ is called {\em stratified} if it maps stratum into stratum, and if 
the restriction of $f$ to each stratum is a submersion.
\end{definition}
\begin{definition}
A stratified map $f:X \to Y,\ X \subset \RM^k$, satisfies the {\em condition $a_f$}
if for any sequence of point $(x_i)$ in a stratum $X_j$ converging to a point $x$ in an
adjacent stratum $X_{j'}$, for which the affine subspaces $\Ker\, df_{|X_j}(x_i)$ converge to a limit $A \subset \RM^k$,
we have the inclusion $\Ker\, df_{|X_{j'}}(x) \subset A$.
\end{definition}
A map germ satisfies the $a_f$-condition if it admits a stratified representative satisfying the $a_f$-condition.
These definitions are due to Thom \cite{Thom_strates} (see also  \cite{Mather_strates}).
The following result is fundamental in the theory of stratified mappings. It is a direct consequence of
 \cite{Gibson_Looijenga}, Chapter II, Theorem 2.6, the proof  being similar to that of  \cite{Gibson_Looijenga}, Chapter II, Theorem
3.2.
\begin{theorem}[\cite{Thom_strates}]
\label{T::Thom}
Let $f:X \to S$ be a stratified mapping satisfying the $a_f$-condition. Let $M$ be a stratum in the Whitney stratification
of $X$.  Any stratified vector field on $M$ tangent to the fibres of $f_{|M}$
extends to a stratified vector field $\theta$ on the strata to which $M$ is adjacent,
tangent to the fibres of the map $f$.
\end{theorem}

\begin{definition}[\cite{Le_oslo}]
A {\em standard representative} $ f:X \to S$ of a flat holomorphic map germ $f:(\CM^n,0) \to (\CM^k,0)$ satisfying the $a_f$ condition
 is a proper stratified representative such that there exists another flat holomorphic representative $ g:Y \to S$
satisfying the following conditions:
\begin{enumerate}
\item $S$ is a closed polydisk containing the origin such that the fibres of $\bar g$ intersect transversally the
boundary of some ball $B_\e$ above $S$ and $X= g^{-1}(S) \cap \dot B_\e$,
\item the fibre $ g^{-1}(0)$ is transverse to the boundary of the balls $B_\e'$ for any $\e' \leq \e$. 
\item the map $g$ satisfies the $a_f$-condition.
\end{enumerate}
\end{definition}
Remark that we abusively denote by the same letter the germ and a standard representative of it.
\subsection{Relative constructibility}
\begin{definition} Consider a stratified map $f:X \to S$ satisfying the $a_f$ condition.
A  sheaf $\Ft$ is called {\em  $f$-constructible} if the following condition holds:
 for each point $x \in X$ there exists a neighbourhood $U$ inside the stratum of $x$
such that
$$\Ft_{|U} \approx f^{-1}(f_{|U})_* \Ft. $$
\end{definition}
If $f$ is the map to a point, a  $f$-constructible sheaf is a  constructible sheaf in the usual sense.
\begin{definition} Consider a stratified holomorphic mapping $f:X \to S$, $X \subset \CM^n$, $S \subset \CM^k$,
satisfying the $a_f$ condition. A complex  $(\Kt^\cdot,\dt) $ of $\Ot_{X}$-coherent
sheaves is called {\em  $f$-constructible} if its cohomology sheaves $\Ht^k(\Kt^\cdot)$ are $f$-constructible and if its differential
is $f^{-1}\Ot_S$ linear.
\end{definition}
The notion of $f$-constructibility extends to germs: given a holomorphic map-germ satisfying the $a_f$-condition
$f:(X,0) \to (S,0)$, a complex $(K^\cdot,\dt) $ of $\Ot_{\CM^{n},0}$-coherent modules is called
{\em  $f$-constructible} if there exists a standard representative $f:X \to S$
and a complex $(\Kt^\cdot,\dt) $ of $f$-constructible $\Ot_{X}$-coherent sheaves such that $K^\cdot$ is the stalk at the origin
of the complex $\Kt^\cdot$.\\
\subsection{The relative de Rham complex}
We give a simple example of relative constructibility: the relative de Rham complex for an isolated
singularity.\\
We consider the relative de Rham complex $\Omega^\cdot_f$ associated to a holomorphic map  $f:X \to S,\ X \subset \CM^n, S \subset \CM^k$
and assume that $S$ is a smooth complex manifold (\cite{Grothendieck_deRham}).\\
For instance if $S = \CM$ then the complex has terms
$\Omega_f^j=\Omega^j_X/\Omega^{j-1}_X \w df$ and the differential is induced by the de Rham differential.
The differential of the complex is obviously $f^{-1}\Ot_{S}$ linear:
$$\pi(d(f\a))=\pi(df \w \a+f d\a)=\pi(fd\a) $$
where $\pi:\Omega_X^\cdot \to \Omega^\cdot_f$ denotes the canonical projection.
A flat holomorphic map-germ $f:(X,0) \to (\CM^k,0)$
defines an {\em isolated singularity} if its special has an isolated singular point at the origin and if it satisfies the
$a_f$-condition.
\begin{proposition}
\label{P::constructible}
The relative de Rham complex $\Omega^\cdot_f$
associated to a flat holomorphic map-germ $f:(X,0) \to (\CM^k,0)$ defining an isolated singularity is $f$-constructible.
\end{proposition}
\begin{proof}
Take a standard representative $f:X \to S$ of the morphism.
Consider the stratification consisting of the smooth points of the morphism (smooth points of the fibres) and its complement.
Stratify the map $f:X \to S$ by refining this stratification.
At a regular point of $f$, the implicit function theorem
shows that the relative de Rham complex is a resolution of the sheaf
$f^{-1}\Ot_S$ (this statement is known as
the relative Poincar\'e lemma).\\
 As the singular set of the fibres is either empty of finite, on each fibre,
the cohomology sheaves restricted to the other strata are constant (any sheaf restricted to a point is constant). 
Therefore the relative de Rham complex is constructible.
\end{proof}
For non isolated singularities the situation if of course much more delicate as the de Rham complex is not always constructible. Additional conditions implying the constructibility have been given in \cite{Barlet, MSaito, vanStraten}.\\
The notions of Theorem \ref{T::Grauert} are now explained.\\

 From Theorem \ref{T::Grauert} and Proposition \ref{P::constructible}, one deduces the following result.
\begin{theorem}
\label{T::Grauert0}
 {For any  flat holomorphic map germ $f:(X,0) \to (\CM^k,0)$ defining
 an isolated singularity, the cohomology spaces $H^p(\Omega_f^\cdot) $ of the relative de Rham complex are
$f^{-1}\Ot_{\CM^k,0}$-coherent modules, for any $p \geq 0$. Moreover, for any  standard representative $f:X \to S$ of the germ $f$, the sheaves
$\RM^\cdot f_*\Omega_f^\cdot$ are $\Ot_S$-coherent and the canonical map $\Omega_f^\cdot(X) \to \Omega_{f,0}^\cdot$
induces an isomorphism $(\RM^p f_*\Omega_f)_0 \approx H^p(\Omega_{f,0}^\cdot) $.}
\end{theorem}
If $X$ is smooth and if the components of $f$ define a complete intersection, these modules are all zero except for $p=0,\dim X-1$
(\cite{Br3,Greuel,Malgrange}).

\section{proof of Theorem \ref{T::Grauert}}
\label{S::proof}
\subsection{Construction of the contraction}
As we saw for the absolute case, the first task is to construct a contraction.
Let $f:X \to S$ be a standard representative, $X =Y \cap \dot B_{\e}$.
\begin{proposition}
\label{P::contraction2}
For any standard representative $f:X \to S$ and any $\e'<\e$,
 the restriction mapping $r:\Kt^\cdot(X) \to \Kt^\cdot(X'),\ X'=f^{-1}(S) \cap \dot B_{\e'}$ is a quasi-isomorphism.
\end{proposition}
\begin{proof}
Let $M \subset Y$ be the stratum of minimal dimension which intersects the ball $B_\e$; and
denote by $\theta_M$ the gradient vector field of the distance function to
$B_\e$ restricted the fibres of $f_{|M}$.\\
By Theorem \ref{T::Thom} the vector field $\theta_M$ extends to a stratified vector field $\theta$
tangent to the fibres of $f$, transversal to the boundary
of $X$ (because of the $a_f$ condition). 
Chose an acyclic covering $U=(U_i)$ of $X$, and take $\e'$ sufficiently close to $\e$, so that
the strata in $U_i'=U_i \cap B_{\e'}$ are homeomorphic to those in $U_i$.\\
As the complex of sheaves is $f$-constructible and $f(U_i)=f(U_i')$,
we have vector space isomorphisms
$$\Ht^q(\Kt^\cdot)(U_i)\approx (f_{|U_i})_*\Ht^q(\Kt^\cdot)f(U_i) \approx \Ht^q(\Kt^\cdot)(U_i') $$ on each small open subset $U_i$. Consequently, the corresponding
hypercohomology spectral sequences  show that the restriction mapping $r:\Kt^\cdot(X) \to \Kt^\cdot(X')$ is a quasi-isomorphism for
any $\e'<\e$.
\end{proof}

\subsection{Fr\'echet modules and nuclear mappings}
\label{SS::nuclear}
Let $A$ be a Fr\'echet algebra, that is a Fr\'echet vector space  for which the multiplication is continuous and associative.
We shall assume that the semi-norms defining the topology also satisfy the inequalities
$p_n(ab) \leq p_n(a)p_n(b),\ \forall n \in \NM$. In the paper of Kiehl-Verdier
this condition is replaced by a more general condition.\\ 
Let us be given two Fr\'echet spaces $E,F$ which are
$A$-modules, that is, modules in the usual sense with the additional condition that the multiplication
mapping is continuous.\\
The topological tensor product $E \hat \otimes_A F$ is defined as the cokernel of the map
$$ (E \hat \otimes_\CM A \hat \otimes_\CM F) \to  (E \hat
\otimes_\CM F),m \otimes a \otimes n \to ma \otimes n-m\otimes an.  $$
A bounded linear mapping 
$u:E \to F$ is called {\em $A$-nuclear}
if it lies in the image of the canonical (Kronecker) embedding $E' \hat \otimes_A
F \to L(E,F)$. Here $E'$ is the strong dual of $E$. These definitions are due to Kiehl-Verdier (\cite{Verdier}).\\
\begin{example}
 Take $A=\Ot_{\CM}(S)$, $E=\Ot_{\CM^2}(D \times S)$, $F=\Ot_{\CM^2}(D' \times S)$ where
$D,D',S$ are disks centred at the origin in $\CM$ such that the respective radius $R,r$ of the disks $D,D' \subset \CM$
are such that $R>r$.\\
The restriction mapping $\rho:\OM_{\CM^2}(D \times S) \to \OM_{\CM^2}(D')$ is $\Ot(S)$-nuclear.
Indeed, let us define the $\Ot(S)$-linear forms $a_n$ 
$$a_n:\Ot_{\CM^2}(D \times S) \to \Ot_{\CM^2}(S),\ f \mapsto \frac{1}{2i\pi}\int_\g \frac{f(z)}{z^{n+1}} dz $$
where $\g$ is a path in $D \setminus D'$ turning counter-clockwise around the origin.
We have
$$\rho=\sum_{n \geq 0} \l^k a_n \otimes \frac{z^n}{\l^n},\ \l=\frac{r}{R}.$$ 
This can also be seen directly using the isomorphism
$$\Ot_{\CM^2}(D \times S) \approx \Ot_{\CM}(D) \otimes \Ot_{\CM}(S).$$
We saw that the restriction mapping $$r: \Ot_{\CM}(D) \to \Ot_{\CM}(D') $$
is $\CM$-nuclear, therefore by tensoring both sides by $ \Ot_{\CM}(S)$, we get an $ \Ot_{\CM}(S)$-nuclear map.
\end{example}
Following Kiehl-Verdier, we say that a map $u:E \to F$ between two $A$-Fr\'echet algebras is {\em $A$-quasinuclear} if there is a commutative diagram
$$\xymatrix{G \ar[rd]^-v \ar[d]^-\pi \\
E \ar[r]^-u & F} $$
where $\pi$ is surjective and $v$ is $A$-nuclear. The reason for considering this notion is explained by
the following proposition which is a direct consequence of Proposition \ref{P::nuclear}.
\begin{proposition}[\cite{Verdier}]
\label{P::quasinuclear}
For any coherent analytic sheaf $\Ft$ in $\CM^n \times \CM^k$
and any Stein neighbourhoods $U,U'$ so that the closure of $U'$ is a compact
subset of $U \subset \CM^n$ and any Stein neighbourhood $S \subset \CM^k$, the restriction mapping $\Ft(U \times \Ot(S)) \to \Ft(U' \times \Ot(S))$ is $\Ot(S)$-quasinuclear.
\end{proposition}
A locally convex vector space $E$ is called {\em nuclear} if any mapping from $E$ to a
Banach space is nuclear. 
\begin{theorem}[\cite{Grothendieck_these,Grothendieck_PTT}]
The Fr\'echet space of holomorphic functions on a polydisk is
a nuclear space.
\end{theorem}
Kiehl and Verdier showed that the following theorem is a consequence of Grothendieck's characterisation of nuclear spaces
(namely the coincidence of the inductive topological tensor product with the projective one).
\begin{theorem}[\cite{Verdier}]
\label{T::nuclear}
For any nuclear Fr\'echet space $E$, the functor $\hat \otimes_\CM E$ is
an exact functor from the category of Fr\'echet spaces to itself.
Moreover, if $E$ is a free module over a
nuclear Fr\'echet algebra\footnote{That $E$ is free means that it is
isomorphic to a product $A \hat \otimes_\CM V$ where $V$ is a nuclear $LF$-space.}$A$,
then the functor $\hat \otimes_A E$ is also exact.
\end{theorem}

\begin{corollary}[\cite{Houzel,Verdier}]
\label{C::produit}
{For any $\Ot_{\CM^{n+k}}$-coherent sheaf $\Ft$, any polycylinders
$S,S' \subset \CM^k$, $S' \subset S$, and any Stein open subset $U \subset \CM^n$,
we have an isomorphism of Fr\'echet modules over the ring 
$\Ot_{\CM^{k}}(S)$
$$  \Ft(U \times S) \hat \otimes \Ot_{\CM^k}(S') \approx
\Ft(U \times S') .$$}
\end{corollary}

\subsection{Comparison of Kiehl-Verdier and Houzel's theorems}
We now formulate a special case of the Houzel theorem.
First, we translate almost literary the following theorems from their original papers the results of Kiehl-Verdier and Houzel.
Most of the notions introduced in both theorems will not be used in this paper.
\begin{theorem}[\cite{Verdier}]
Given a nuclear chain $(A_t,\p_t)$ of algebras and a quasinuclear quasi-isomorphism between  complexes bounded from above
$$f^\cdot:N^\cdot \to M^\cdot$$
such that $M^\cdot$ is transversal to $A_t$ for all $t$'s. $M^n$ and $N^n$ are nuclear for all $n$. $\p_{t't}$ satisfies condition $B$
for all $t'<t$. Then there exists a complex of finite type free $A_\beta$-modules, bounded from above, and a quasi-isomorphism
$$F^\cdot \to M^\cdot.  $$ 
\end{theorem}
The following theorem due to Houzel shows that most assumptions of the
Kiehl-Verdier theorem are superfluous.
 \begin{theorem}[\cite{Houzel}]
\label{T::Houzel}Let $A$ be a multiplicatively convex, complete bornological algebra and let $(M_i^\cdot)$ be a sequence
of complexes of complete bornological $A$ modules; let be given for $1 \leq i \leq r$, a homomorphism of complexes
$u_i:M_{i-1}^\cdot \to M_i^\cdot$ ($A$-linear and bounded), and integers $a,b \in \ZM$, with $a \leq b$.\\
We make the following assumptions
\begin{enumerate}
\item for all $i$, $M_i^n$ satisfies the homomorphism property for $n \geq a$ and is zero for $n \geq b$,
\item for $1 \leq i \leq r$, $u_i$ is a quasi-isomorphism (see \cite{Illusie}) and is $A$-nuclear in degree $\geq a$,
\item $r \geq b-a+1$.
\end{enumerate}
Then the complexes $M_i^\cdot$ are $a$-pseudo-coherent (see \cite{Illusie}).
\end{theorem}
We now deduce the Kiehl-Verdier from the Houzel theorem by reducing the above theorem to some special case.\\
In our case, we take $r=3$, so that $i=0,1,2$, and $a=0$, $b=2$, so that $M_2$ is zero, The complex denoted $N^\cdot,M^\cdot$ in the
Kiehl-Verdier theorem are denote  $M_0^\cdot,M_1^\cdot$ in the formulation given by Houzel.\\
The Houzel theorem implies the following result\footnote{We slightly changed the formulation by replacing nuclear by quasinuclear but this does not change anything in Houzel's proof, for obvious reasons.}. 
 \begin{theorem}
\label{T::Houzelbis}Let $A$ be a Fr\'echet algebra and $M^\cdot, N^\cdot$ two complexes of Fr\'echet spaces
for which there exists a quasinuclear quasi-isomorphism $u:N^\cdot \to M^\cdot$.
Then, there exists a complex of finite type $A$-modules quasi-isomorphic to $M^\cdot$ and $N^\cdot$.
\end{theorem}
It is a remarkable fact that the proof of this result relies essentially on a generalisation of Lemma \ref{L::Houzel} to Fr\'echet modules.
\subsection{Proof of Theorem \ref{T::Grauert}}
\begin{lemma} The restriction mapping $r:\Kt^\cdot(X) \to \Kt^\cdot(X')$
is a $\Ot_{S}(S)$-quasinuclear quasi-isomorphism.
\end{lemma}
\begin{proof}
Recall that $X$ is the intersection of a Stein open neighbourhood with some open ball $B_\e$.
Define the complex of sheaves $\t \Kt^\cdot$ in $B_\e \times S$ by the presheaf
$$\t \Kt^\cdot(U \times V)= \Kt^\cdot(U \cap f^{-1}(V)).$$
Both complexes are isomorphic as complexes of Fr\'echet sheaves. Moreover, the restriction mapping 
$\t \Kt^\cdot(B_\e \times S) \to \t \Kt(B_{\e'} \times S)$ is $\Ot_S(S)$-quasinuclear (Proposition \ref{P::quasinuclear}). This proves the lemma.
\end{proof}
This lemma shows that Theorem \ref{T::Houzelbis} applies, therefore there exists a complex $\Lt^\cdot$ of free coherent
$\Ot_{S}$-sheaves so that $\Lt^\cdot(S)$ is quasi-isomorphic to $\Kt^\cdot(X)$.
\begin{lemma} The sheaf complexes $\Lt^\cdot$, $ f_* \Kt^\cdot_{|X}$ are quasi-isomorphic
\end{lemma}
\begin{proof}
A mapping $u:M^\cdot \to L^\cdot$ of complexes induces a quasi-isomorphism between
two complexes if and only if its mapping cylinder $C^\cdot(u)$ is exact.\\
We apply this fact to the mapping cylinder of the quasi-isomorphism 
$$u:\Lt^\cdot(S) \to \Kt^\cdot(X).$$ 
As the functor $\hat \otimes \Ot_S(P)$ is exact for any polydisk $P \subset S$ (Theorem \ref{T::nuclear}),
the complex $C^\cdot(u) \hat \otimes \Ot_{S}(P)$ is also exact.\\
Using Corollary \ref{C::produit}, we get that
 $C^\cdot(u) \hat \otimes  \Ot_{S}(P)$ is the mapping cylinder of
$u':\Lt^\cdot(P) \to \Kt^\cdot(X \cap f^{-1}(P))$.
Therefore, the complexes of sheaves $\Lt^\cdot$ and $ f_* \Kt^\cdot_{|X}$ are
quasi-isomorphic. This proves the lemma.
\end{proof}
I assert that the complex $K^\cdot=\Kt^\cdot_0$ is quasi-isomorphic to the stalk of the complex $\Lt^\cdot$ at the origin.\\
Let $(B_{\e_n})$ be a fundamental sequence of neighbourhoods of the origin in $\CM^n$,
so that their intersection with the special fibre of $f$
is transverse. As the map $f$ satisfies the $a_f$-condition, we may find a fundamental sequence $(S_n)$ of neighbourhoods of the origin
in $\CM^k$ so that the fibres of $f$ intersect $B_{\e_n}$ transversally above $S_n$.\\
Put $X_n=f^{-1}(S_n)$, we have the isomorphism
$$\Lt_{|S_n}^\cdot \approx f_* \Kt^\cdot_{|X_n} \approx  f_* \Kt^\cdot_{|X_n \cap B_{\e_n}}. $$
The first isomorphism is a consequence of the previous lemma and the second one follows from the fact that the
contraction is a quasi-isomorphism (Proposition \ref{P::contraction2}).\\
In the limit $n \to \infty$, we get that the complex $K^\cdot=\Kt^\cdot_0$ is quasi-isomorphic to the complex $\Lt^\cdot_0$.
This concludes the proof of the theorem.
\section{Finiteness theorem for coherent ind-analytic complexes}
\subsection{The sheaves $\Ot_{X|Y}$}
Let $i:X \to Y$ be the inclusion of a complex analytic manifold $X$ into another complex analytic manifold $Y$.
The sheaf $i^{-1}\Ot_{Y}$ is denoted by $\Ot_{X|Y}$. If $Y$ is of the form $X \times T$, we denote simply by $\Ot_{X|X \times T}$
the sheaf obtained from the inclusion of $X \times \{ 0 \}$ in $X \times T$. These sheaves are frequently considered in
microlocal analysis (\cite{SKK}).\\
The stalk of  the sheaf $\Ot_{X|X \times T}$ at a point $x_0$ is the space of germs of holomorphic functions in $X \times T$
at the point $x=x_0,\ t=0$.\\
It follows from Cartan's theorem A that the ring  $\Ot_{X|X \times T}$ is coherent, that is, the kernel of any map
$$\Ot^k_{X|X \times T} \to \Ot_{X|X \times T}$$
is finitely generated. Following Serre \cite{Serre_FAC}, we say that a sheaf $\Ft $ on a space $X$ is {\em $\Ot_{X|X \times T}$-coherent}
or that it is a {\em coherent ind-analytic sheaf}, if it is the cokernel of a morphism of $\Ot_{X|X \times T}$-modules:
$$\Ot_{X|X \times T}^{p} \to \Ot_{X|X \times T}^{n} \to \Ft \to 0.$$ 
The notion of $f$-constructibility extends trivially to complexes of
$\Ot_{X|X \times T}$-coherent sheaves and to their stalks.
\subsection{Topological properties of the sheaves  $\Ot_{X|Y}$}
\label{SS::Dieudonne_Schwartz}
The sheaves  $\Ot_{X|Y}$ are inductive limits of Fr\'echet spaces,  called {\em $LF$-spaces}. A notion that we will now recall.
Consider a set of continuous
morphisms of Fr\'echet spaces to a fixed vector space
$u_i:E_i \to E$, $i \in \Omega$ such that $\bigcup_{i \in \Omega} u_i(E_i)=E$. The {\em inductive limit topology} $T$ of the vector space $E$ is defined by
$$ U \in T \iff \forall i \in \Omega,\ u_i^{-1}(U) {\rm \ is\ open \ in\ } E_i.$$
Such topological vector spaces are locally convex and bornological, i.e., bounded linear mappings coincide with continuous ones.\\
In case the $E_i$'s form an increasing sequence of closed vector subspaces in $E$ and the $u_i$'s are the inclusions, the resulting $LF$-space is complete and satisfies the Banach theorem, that is,
bounded surjective linear mappings are open (\cite{Dieudonne_Schwartz}).
In the general case, $LF$-spaces are not always complete but if they are, then
they satisfy the Banach theorem provided that the set $\Omega$ is countable (\cite{Koethe_cexample,Koethe_Banach}).\\
The space $\Ot_{\CM^n}(A)=\underrightarrow{\lim}\Ot_{\CM^n}(U),\ A \subset U$ of holomorphic
functions restricted to a compact subset $A \subset \CM^n$ is a $LF$-space, for instance the germs at a point
of holomorphic functions in $\CM^n$ form an $LF$-space. This space is complete, Montel and reflexive (\cite{Grothendieck_EVT}, Chapter 3, Proposition 5, Corollary 2 and Example b).\\
In particular, for any compact neighbourhood $A \subset X$, the algebra $\Ot_{X|Y}(A)$ satisfies the assumptions of the Houzel theorem.\\
As inductive limits commute with topological tensor products, the topological vector space  $\Ot_{X|Y}(A)$ is nuclear.
\subsection{Local finiteness theorem}
As the Houzel theorem applies for coherent ind-analytic sheaves, the proof of the following theorem is exactly analogous to that of the analytic case (from now on, we do no longer consider direct image sheaves as the conclusions are in all cases similar). 
\begin{theorem}
\label{T::Grauert2}
 {Let $f:(\CM^n,0) \to (\CM^k,0)$ be a holomorphic map germ satisfying the $a_f$-condition.
The cohomology spaces $H^p(K^\cdot) $ associated to a complex of  $f$-constructible
$\OM_{\CM^n|\CM^{n+N},0}$-coherent modules are $f^{-1}\Ot_{\CM^k|\CM^{k+N},0}$-coherent modules, for any $p \geq 0$.}
\end{theorem}

Let us denote by $\Mk$ the maximal ideal of the local ring $\Ot_{\CM^k,0}$. We identify $\CM^k$ to
$\CM^k \times \{ 0 \} \subset \CM^k \times \CM^n, $
From the above theorem, one deduces the following algebraic form of the division theorem due to Houzel and Serre.
\begin{theorem}[(\cite{Houzel_division}]
Any $\Ot_{\CM^{k+n},0}$-module  $M$  of finite type such that $M/\Mk M$ is a finite dimensional
$\CM$-vector space is itself an $\Ot_{\CM^k,0}$-module of finite type.
\end{theorem}
\begin{proof}
 Consider the complex obtained from
a projective resolution of $M$:
$$K^\cdot:\xymatrix{ \dots \ar[r]^-{\dt_2} & \Ot_{\CM^{k+n},0}^{n_1} \ar[r]^-{\dt_1}& \Ot_{\CM^{k+n},0}^{n_0} \ar[r] & 0 },\ H^0(K^\cdot)=\Ot_{\CM^{k+n},0}^{n_0}/{\rm Im}\, \dt_1 \approx M.$$
This complex gives a complex $\Kt^\cdot$ of $\Ot_{X \times T}$-sheaves which support is an analytic variety $V \subset X \times T$.
Here $X \subset \CM^n$ and $T \subset \CM^k$ are small neighbourhoods of the origin in $\CM^n$ and $\CM^k$.\\
The dimension of $M/\Mk M$ is equal to the intersection multiplicity of $V$ with the affine subspace
$X \times \{ 0 \}=\{ (x,t), t=0 \}$. Assume that $X$ and $T$ are small enough so that the variety $V$ intersects
the affine space $X \times \{ 0 \}$ only at the origin.
Then, the restriction of the complex $\Kt^\cdot$ to $X \times \{0 \}$ defines a $\Ot_{X|X \times T,0} $-coherent sheaf complex in a neighbourhood of the origin $X \subset \CM^n$.
The cohomology of the complex is supported at the origin, it is therefore constructible.
Theorem \ref{T::Grauert} implies 
that $H^p(K^\cdot)=M$ is a finite type module over the ring $\Ot_{\CM^k,0}$ (here $f$ is the mapping to a point).
\end{proof}
\section{Finiteness theorem for non-commutative sheaves}


\subsection{Statement of the theorem}
We now give a generalisation of Theorem \ref{T::Grauert} which applies to non-commutative
algebras such as algebras of pseudo-differential operators.\\
Given a holomorphic function $f:X \to S$ between open subsets $X\subset \CM^n,\ S \subset \CM^k$, we define a category $\Ht(f)$.\\
The objects of the category $\Ht(f)$ are triples
$(\Kt^\cdot,\At,\Bt)$ where  $\At$ and $\Bt$ are sheaves of $LF$-algebras in the closed neighbourhoods $X$
and $S$, the complex $\Kt^\cdot$ is an $f$-constructible sheaf complex of coherent $\At$-modules in $X$, the sheaves
$f^{-1} \Bt \subset \At$ are sheaves of subalgebras  and the differential of the complex is $f^{-1}\Bt$-linear.\\
A morphism between two objects is a continuous linear mapping between the complexes and the algebras without regarding the module
and ring structures.\\
We say that a triple
is {\em $f$-holomorphic},  if it consists of a
complex of $\Ot_{X|X \times T}$-coherent sheaves with a
$f^{-1}\Ot_{S|S \times T}$ linear differential, for some complex reduced analytic space $T$.\\
Holomorphic triples satisfy the assumption of the Houzel theorem and this property is preserved under isomorphism, therefore we have the following result.
\begin{theorem}
\label{T::Grauert3}
Assume that the function
$f:X \to S$ satisfies the $a_f$-condition. For all triples
$(\Kt^\cdot,\At,\Bt) \in \Ht(f)$ isomorphic to an $f$-holomorphic triple, the cohomology sheaves $\Ht^p(\Kt^\cdot)$
are sheaves of coherent $f^{-1}\Bt$-modules, for any $p \geq 0$.
\end{theorem} 

\subsection{The Boutet de Monvel theorem}
In \cite{SKK}, the authors proved that a division theorem for
pseudo-differential operators.
In \cite{BdM}, Boutet de Monvel showed that this theorem admits an algebraic version
(or rather a generalisation) similar to the Houzel-Serre formulation of the
division theorem  (see also \cite{Pham_livre}, Chapter 3). We give a simple proof of this theorem based on the absolute case of Theorem \ref{T::Grauert3}.\\
Denote by $\Et(0)$ the sheaf of analytic pseudo-differential operator in
$T^*\CM^n \approx \CM^{2n}=\{(q,p) \}$ of order $0$.
Let $\Et'(0)$ be a subsheaf of operators
which depend only on some of the variables, say
$q_1,\dots,q_j,p_1,\dots,p_k$.
We denote by  $\Et_0(0),\Et'_0(0)$ the stalks of the sheaves  $\Et(0),\Et'(0)$
at the point $x_0 \in \CM^{2n}$ with coordinates $q_1=\dots=q_n=0$, $p_1=1,p_2=0,\dots,p_n=0$.
\begin{theorem}[\cite{BdM,Pham_livre}]
For any coherent left $\Et_0(0)$-module  $M$ the following assertions are
equivalent
\begin{enumerate}
\item the $\Ot_{\CM^{j+k-1},0}$-module $M/\d_{q_1}^{-1}M$ is of finite type,
\item  the  $\Et'_0(0)$-left module $M$ is of finite type
\end{enumerate}
\end{theorem}
\begin{proof}
Denote by $X \approx \CM^{n-j} \times \CM^{k-1}$ the vector subspace of $\CM^{2n}$ defined by the equations
$$q_1=\dots=q_j=p_2=\dots=p_k=0.$$
The module $M$ is the stalk at the point $x_0=(0,\dots,0,1,0,\dots,0)$ of a sheaf $\Mt$ of $\Et_0(0)$-modules in
$T^*\CM^n \approx \CM^{2n}$.\\
Consider the complex given by a resolution of $\Mt$
$$\Kt^\cdot:\xymatrix{\cdots \ar[r]^-{\dt_2}& \Et(0)^{n_1} \ar^-{\dt_1}[r]& \Et(0)^{n_0} \ar[r] & 0 },\
\Ht^0(\Kt^\cdot)= \Et(0)^{n_0}/{\rm Im}\, \dt_1 \approx \Mt.$$
The  support of $\Mt$ coincides with that of $\Mt/\d_{q_1}^{-1}\Mt$, it is therefore
an analytic subvariety $V \subset \CM^{2n}$ \cite{SKK} (see also \cite{Pham_livre}, Proposition 4.2.0).\\
As the $\Ot_{\CM^{j+k-1},0}$-module $M/\d_{q_1}^{-1}M$ is of finite type,
in a sufficiently small neighbourhood $U$ of the origin, the variety $V$ intersects the vector space
$X$ only at the origin. Therefore the complex $\Kt^\cdot$ restricted to $X \cap U$ is constructible
for the Whitney stratification of $X \cap U$ consisting of the origin and its complement.
The restriction of the sheaf $\Et(0)$ to the complement of the zero section in $T^*\CM^n \approx \CM^{2n}$
is a sheaf of Fr\'echet algebras (\cite{Boutet_Kree}); consequently Theorem \ref{T::Grauert3} applies. This concludes
the proof of the theorem.
\end{proof}


\bibliographystyle{amsplain}
\bibliography{master}

\providecommand{\bysame}{\leavevmode\hbox to3em{\hrulefill}\thinspace}
\providecommand{\MR}{\relax\ifhmode\unskip\space\fi MR }
\providecommand{\MRhref}[2]{%
  \href{http://www.ams.org/mathscinet-getitem?mr=#1}{#2}
}
\providecommand{\href}[2]{#2}
\begin{thebibliography}{10}

\bibitem{Barlet}
D.~Barlet and M.~Saito, \emph{{Brieskorn modules and Gauss-Manin systems for
  non-isolated hypersurfaces singularities}}, ArXiv:math.CV/0411406, 17 p.

\bibitem{Bourbaki_evt}
N.~Bourbaki, \emph{Espaces vectoriels topologiques}, Hermann, 1966.

\bibitem{BdM}
L.~Boutet~de Monvel, \emph{Op\'erateurs pseudo-diff\'erentiels analytiques},
  S\'eminaire op\'erateurs diff\'erentiels et pseudo-diff\'erentiels, Grenoble.

\bibitem{Boutet_Kree}
L.~Boutet~de Monvel and P.~Kr\'ee, \emph{{Pseudo differential operators and
  Gevrey classes}}, Annales de l'Institut Fourier \textbf{17} (1967), no.~1,
  295--323.

\bibitem{Br3}
E.~Brieskorn, \emph{{Die Monodromie der isolierten Singularit\"aten von
  Hyperfl\"achen}}, Manuscr. Math. \textbf{2} (1970), 103--161.

\bibitem{Buchweitz}
R.O. Buchweitz and G.M. Greuel, \emph{The {Milnor} number and deformations of
  complex curve singularities}, Invent. Math. \textbf{58} (1980), 241--281.

\bibitem{Cartan_Serre}
H.~Cartan and J.P. Serre, \emph{Un th\'eor\`eme de finitude concernant les
  vari\'et\'es analytiques compactes}, Comptes Rendus \`a l'Acad\'emie des
  Sciences \textbf{237} (1953), 128--130.

\bibitem{Dieudonne_Schwartz}
J.~Dieudonn\'e and L.~Schwartz, \emph{{La dualit\'e dans les espaces $(\mathcal
  F)$ et $(\mathcal{L}\mathcal{F})$}}, Annales de l'institut Fourier \textbf{1}
  (1949), 61--101.

\bibitem{Douady}
A.~Douady, \emph{{Le th\'eor\`eme des images directes de Grauert (d'apr\`es
  Kiehl-Verdier)}}, Ast\'erisque \textbf{16} (1974), 49--62.

\bibitem{DouadyR.}
R.~Douady, \emph{{Produits tensoriels topologiques et espaces nucl\'eaires}},
  Ast\'erisque \textbf{16} (1974), 1--25.

\bibitem{Forster_Knorr}
O.~Forster and K.~Knorr, \emph{{Ein Beweis des Grauertschen Bildgarbensatzes
  nach Ideen von B. Malgrange}}, Manuscripta Math. \textbf{5} (1971), 19--44.

\bibitem{Gibson_Looijenga}
C.G. Gibson, K.~Wirthm\"uller, A.A. du~Plessis, and E.J.N. Looijenga,
  \emph{Topological stability of smooth mappings}, Lecture Notes in
  Mathematics, vol. 552, Springer, 1976.

\bibitem{Greuel}
G.-M. Greuel, \emph{{Der Gauss-Manin Zusammenhang isolierter Singularit\"aten
  von vollst\"andigen Durchschnitten}}, Math. Ann. \textbf{214} (1975),
  235--266.

\bibitem{Grothendieck_these}
A.~Grothendieck, \emph{{R\'esum\'e des r\'esultats essentiels dans la th\'eorie
  des produits tensoriels topologiques et des espaces nucl\'eaires}}, Annales
  de l'institut Fourier (1952), 73--112.

\bibitem{Grothendieck_PTT}
\bysame, \emph{Produits tensoriels topologiques et espaces nucl\'eaires}, Mem.
  of the Am. Math. Soc. \textbf{16} (1955).

\bibitem{Grothendieck_deRham}
\bysame, \emph{{On the de Rham cohomology of algebraic varieties}},
  Publications Math\'ematiques de l'IH\'ES \textbf{29} (1966), 95--103.

\bibitem{Grothendieck_EVT}
\bysame, \emph{Topological vector spaces}, Gordon and Breach, 1973, 245 p.,
  English Translation: Espaces vectoriels topologiques, S\~ao Paulo 1954.

\bibitem{Houzel_division}
C.~Houzel, \emph{{G\'eom\'etrie analytique locale, I}}, vol.~13, S\'eminaire
  Henri Cartan, no.~2, ch.~Expos\'e No. 18, 1960, 12 p.

\bibitem{Houzel}
\bysame, \emph{Espaces analytiques relatifs et th\'eor\`eme de finitude}, Math.
  Annalen \textbf{205} (1973), 13--54.

\bibitem{Houzel_Schapira}
C.~Houzel and P.~Schapira, \emph{Images directes de modules diff\'erentiels},
  C. R. Acad. Sci. Paris \textbf{298} (1984), no.~18, 461--464.

\bibitem{Hubbard}
J.~Hubbard, \emph{{Transversalit\'e}}, Ast\'erisque \textbf{16} (1974), 33--48.

\bibitem{Illusie}
L.~Illusie, \emph{G\'en\'eralit\'es sur les conditions de finitude dans les
  cat\'egories d\'eriv\'ees}, SGA VI, Lecture Notes in Mathematics, vol. 225,
  Springer, 1971, pp.~78--159.

\bibitem{Verdier}
R.~Kiehl and J.L. Verdier, \emph{{Ein Einfacher Beweis des Koh\"arenzsatzes von
  Grauert}}, Math. Annalen \textbf{195} (1971), 24--50.

\bibitem{Koethe_cexample}
G.~K\"othe, \emph{{\"Uber die Vollst\"andigkeit einer Klasse lokalkonvexer
  R\"aume}}, Mathematische Zeitschrift \textbf{52} (1950), 627--630.

\bibitem{Koethe_Banach}
\bysame, \emph{{\"Uber zwei S\"atze von Banach}}, Mathematische Zeitschrift
  \textbf{53} (1950), 203--209.

\bibitem{Le_oslo}
D.T. L\^e, \emph{Some remarks on relative monodromy}, Real and complex
  singularities, Oslo 1976 (P.~Holm, ed.), Sijthoff and Noordhoff, 1977,
  pp.~397--404.

\bibitem{Levy}
R.~Levy, \emph{{A new proof of the Grauert direct image theorem}}, Proceedings
  of the American Mathematical Society \textbf{99} (1987), no.~3, 535--542.

\bibitem{Lojasiewicz}
S.~Lojasiewicz, \emph{{Stratification des ensembles analytiques avec les
  propri\'et\'es (A) et (B) de Whitney.}}, Fonctions analytiques de plusieurs
  variables et analyse complexe (Colloq. Internat. C.N.R.S. No. 208, Paris,
  1972), vol.~1, Gauthier-Villars, Paris, 1974.

\bibitem{Malgrange}
B.~Malgrange, \emph{Int\'egrales asymptotiques et monodromie}, Ann. scient.
  \'Ecole Norm. Sup. \textbf{7} (1974), no.~4, 405--430.

\bibitem{Mather_strates}
J.~Mather, \emph{Stratifications and mappings}, Dynamical systems (Proc.
  Sympos., Univ. Bahia, Salvador, 1971), Academic Press, 1973, pp.~195--232.

\bibitem{Pham_livre}
F.~Pham, \emph{{Singularit\'es des syst\`emes diff\'erentiels de Gauss-Manin}},
  Birkh\"auser, 1979, 339 pp.

\bibitem{MSaito}
M.~Saito, \emph{{Gauss-Manin connection for non-isolated hypersurface
  singularities}}, Master thesis, Univ. Tokyo, 1980.

\bibitem{SKK}
M.~Sato, T.~Kawai, and M.~Kashiwara, \emph{Microfunctions and
  pseudodifferential equations}, Hyperfunctions and pseudodifferential
  equation, Lecture Notes in Mathematics, no. 287, Springer, 1971.

\bibitem{Schneiders}
J.-P. Schneiders, \emph{{A Coherence Criterion for Fr\'echet Modules}},
  Ast\'erique \textbf{224} (1994), 99--113, Part I.

\bibitem{Schwartz}
L.~Schwartz, \emph{Homomorphismes et applications compl\`etement continues},
  Comptes Rendus \`a l'Acad\'emie des Sciences \textbf{236} (1953), 2472--2473.

\bibitem{Sebastiani}
M.~Sebastiani, \emph{{Preuve d'une conjecture de Brieskorn}}, Man. Math.
  \textbf{2} (1970), 301--308.

\bibitem{Serre_FAC}
J.-P. Serre, \emph{Faisceaux alg\'ebriques coh\'erents}, Ann. of Math.
  \textbf{61} (1955), 197--278.

\bibitem{Teissier_Oslo}
B.~Teissier, \emph{The hunting of invariants in the geometry of discriminants},
  Real and complex singularities (Proc. Ninth Nordic Summer School/NAVF Sympos.
  Math., Oslo, 1976), Sijthoff and Noordhoff, Alphen aan den Rijn, 1977,
  pp.~565--678.

\bibitem{Te2}
\bysame, \emph{Vari\'et\'es polaires {II}}, Algebraic Geometry Proc. La
  M\'abida, Lect. Notes in Math. n.961, Springer-Verlag, 1981, pp.~314--495.

\bibitem{Thom_strates}
R.~Thom, \emph{Ensembles et morphismes stratifi\'es}, Bull. Amer. Math. Soc.
  \textbf{75} (1969), 240--284.

\bibitem{vanStraten}
D.~van Straten, \emph{{On the Betti numbers of the Milnor fibre of a certain
  class of hypersurfaces singularities}}, Singularities, representation of
  algebras and vector bundles (G.-M. Greuel and G.~Trautmann, eds.), Lecture
  Notes in Math., no. 1273, Springer-Verlag, 1987, pp.~203--220.

\bibitem{Wh}
H.~Whitney, \emph{Tangents to an analytic variety}, Ann. of Math. \textbf{81}
  (1964), 496--549.

\end{thebibliography}
\end{document}